\documentclass[12pt]{article}
\usepackage{amssymb}
\usepackage{amsmath}
\usepackage{indentfirst,latexsym}
\usepackage{mathrsfs}
\usepackage{graphicx}
\usepackage{fancyhdr}
\date{}
\begin{document}

\title{Two permutation classes enumerated by the central binomial coefficients}
\author{Marilena Barnabei, Flavio Bonetti,\\ \ \\
 Dipartimento di Matematica, Universit\`a di Bologna\\ P.zza di Porta San Donato 5, 40126 Bologna, Italy\\ \ \\ Matteo Silimbani,\\ \ \\ LaBRI - 351, Cours de la Lib\'eration, 33405 Talence, France} \maketitle

\noindent {\bf Abstract.} We define a map between the set of permutations that avoid
either the four patterns $3214,3241,4213,4231$ or $3124,3142,4123,4132$, and the set of Dyck prefixes. This map, when restricted to either of the two classes, turns out to be a bijection that allows us to determine some
notable features of these permutations, such as the distribution of the statistics ``number of ascents'', ``number of left-to-right maxima'',
``first element'', and ``position of the maximum element''.


\section{Introduction}

\noindent A permutation $\sigma$ is said to \emph{contain} a permutation
$\tau$ if there exists a subsequence of $\sigma$ that has the same
relative order as $\tau$, and in this case $\tau$ is said to be a
\emph{pattern} of  $\sigma$. Otherwise, $\sigma$ is said to
\emph{avoid the pattern} $\tau$.

\noindent A \emph{class} of permutations is a downset in the
permutation pattern order defined above. Every class can be
defined by its basis $B$, namely, the set of minimal permutations
that are not contained in it; the class is denoted by $Av(B)$. We denote by
$S_n(B)$ the set $Av(B)\cap S_n$.\newline

\noindent The classes of permutations that avoid one or more patterns of length $3$ have been exhaustively studied
since the seminal paper \cite{simisch}. In many cases, the properties of these permutations have been determined by establishing suited bijections with lattice paths (see a survey in \cite{ck}).

\noindent The case of patterns of length $4$ still seems to be incomplete, even with regard to the mere enumeration in the case of multiple avoidance. In his thesis \cite{guib}, O. Guibert deals with some enumerative problems concerning multiple avoidance. In particular, Theorem 4.6 exhibits $12$ different classes of permutations avoiding four patterns of length $4$, each one enumerated by
the sequence of central binomial coefficients.\\

\noindent On the other hand, it is well known that the central binomial coefficient ${2n\choose n}$ enumerates the set of Dyck prefixes of length $2n$, namely, lattice paths in the integer
lattice $\mathbb{N}\times\mathbb{N}$ starting from the origin,
consisting of up steps $U=(1,1)$ and down steps $D=(1,-1)$, and
never passing below the $x$-axis.\\

\noindent In this paper we consider two of Guibert's classes, namely, $Av(T_1)=$\\ $Av(3214,3241,4213,4231)$ and $Av(T_2)=Av(3124,3142,4123,4132)$, and  introduce a map $\Phi$ between the union of these classes and the set of Dyck prefixes. The restrictions of this map to the sets $Av(T_1)$ and $Av(T_2)$ turn out to be bijections.
The key tool in determining this map is Theorem \ref{pincopallo}, which describes the structure of permutations of both classes. If we consider the decomposition $\sigma=\alpha\,n\,\beta$, where $n$ is the maximum symbol in $\sigma$ and $\alpha$ and $\beta$ are possibly empty words, then if the permutation $\sigma$ avoids the four patterns in $T_1$, the prefix $\alpha$ avoids $321$ and the suffix $\beta$ avoids both $231$ and $213$, while if $\sigma$ avoids $T_2$, the prefix $\alpha$ avoids $312$ and the suffix $\beta$ avoids both $123$ and $132$. In both cases, the lattice path $\Phi(\sigma)$ is obtained by first associating with $\alpha$ a Dyck prefix by a procedure similar to the one used by Krattenthaler in \cite{kratt}, and then appending to this prefix a sequence of up steps and down steps that depends on the suffix $\beta$.\\

\noindent The map $\Phi$ allows us to relate some properties of a permutation $\sigma$ with some particular features of the corresponding Dyck prefix $P$. For example, the Dyck prefix $P$ does not touch the $x$-axis (except for the origin) whenever $\sigma$ is connected, while $P$ is a Dyck path whenever $\sigma$ ends with the maximum symbol. Moreover, if $\sigma\in Av(T_1)$, the $y$-coordinate of the last point of $P$ gives information about the maximum length of a decreasing subsequence in $\sigma$.

\noindent The map $\Phi$ allows us also to prove that each one of the three statistics ``number left-to-right maxima'', ``position of the maximum element'', and ``first element'' are equidistributed over the two classes, and we  determine the generating function of these statistics.
Finally, in the last section we determine the distribution of the statistic ``number of ascents'', which is not equidistributed over the two classes.

\section{The classes $Av(3214,3241,4213,4231)$ and \\ $Av(3124,3142,4123,4132)$}

\noindent In this paper we are interested in the two classes $Av(T_1)$ and $Av(T_2)$, where $T_1=\{3214,3241,4213,4231\}$ and $T_2=\{3124,3142,4123,4132\}$. First of all, we characterize the permutations in these classes by means of their
left-to-right-maximum decomposition.\newline

\noindent Recall that a permutation $\sigma$ has a
\emph{left-to-right maximum} at position $i$ if
$\sigma(i)\geq\sigma(j)$ for every $j\leq i$, and that every
permutation $\sigma$ can be decomposed as
$$\sigma=M_1\,w_1\,M_2\,w_2\,\ldots\,M_k\,w_k,$$
where $M_1,\ldots,M_k$ are the left-to-right maxima of $\sigma$
and $w_1,\ldots,w_k$ are (possibly empty) words. The following
characterization of permutations in $Av(T_1)$ and $Av(T_2)$ is easily deduced:

\newtheorem{inizio}{Theorem}
\begin{inizio}\label{pincopallo}
Let $\sigma\in S_n$. Consider the decomposition
$$\sigma=M_1\,w_1\,M_2\,w_2\,\ldots\,M_k\,w_k,$$ where $M_k=n$.
Denote by $l_i$ the length of the word $w_i$. Then:
\begin{itemize}
\item[a.] $\sigma$
belongs to $Av(T_1)$ if and only if: \begin{itemize} \item the
renormalization of $w_k$ is a permutation in $Av(231,213)$, and
\item if $k>1$, the juxtaposition of the words $w_1,\ldots,w_{k-1}$
consists of the smallest $l_1+\cdots+l_{k-1}$ symbols in the set
$[n-1]\setminus\{M_1,\ldots,M_{k-1}\}$, listed in increasing
order. In particular, the permutation\\
$\alpha=M_1\,w_1\,M_2\,w_2\,\ldots\,M_{k-1}\,w_{k-1}$, after
renormalization, avoids $321$.
\end{itemize}
\item[b.] $\sigma$
belongs to $Av(T_2)$ if and only if: \begin{itemize} \item the
renormalization of $w_k$ is a permutation in $Av(123,132)$, and
\item if $k>1$, every word $w_i$, $i\leq k-1$, consists of the $l_i+1$ greatest unused symbols among those that are less than $M_i$, listed in decreasing order. In particular,
$\alpha=M_1\,w_1\,M_2\,w_2\,\ldots\,M_{k-1}\,w_{k-1}$, after
renormalization, avoids $312$.
\end{itemize}
\end{itemize}\end{inizio}
\begin{flushright}
\vspace{-.5cm}$\diamond$
\end{flushright}

\noindent This result implies that a permutation $\sigma\in S_n(T_1)$ can be decomposed into
$$\sigma=\alpha\, n \, \beta,$$
where $\alpha$ avoids $321$ and $\beta$ avoids both $231$ and $213$.
Similarly, a permutation $\sigma'\in S_n(T_2)$  can be decomposed into
$$\sigma'=\alpha'\, n \, \beta',$$
where $\alpha'$ avoids $312$ and $\beta'$ avoids both $123$ and $132$. We submit that, in both cases, this property is not
a characterization, since the permutations $\alpha$ and $\alpha'$ can not be chosen arbitrarily, according to Theorem \ref{pincopallo}:
for example, the permutation $3\,2\,4\,1$, that belongs to $T_1$, has precisely the described structure.\\

\noindent Recall that a permutation
$\tau=x_1\,x_2\,\ldots\,x_j$ belongs to $Av(231,213)$ whenever,
for every $i\leq j$, the integer $x_i$ is either the minimum or
the maximum of the set $\{x_i,x_{i+1},\ldots,x_j\}$. Analogously, $\tau$ belongs to $Av(123,132)$
whenever,
for every $i\leq j$, the integer $x_i$ is either
the greatest or the second greatest element of the set $\{x_i,x_{i+1},\ldots,x_j\}$ (see, e.g.,
\cite{simisch}).\\

\noindent For example, if we consider the permutation in $S_{10}(T_1)$
$$\sigma=4\,1\,2\,6\,7\,3\,10\,5\,9\,8,$$ we have
$$\alpha=4\,1\,2\,6\,7\,3,$$
with
$$M_1=4\quad M_2=6\quad M_3=7\quad M_4=10,$$
and
$$\beta=5\,9\,8.$$
Analogously, if we consider the permutation in $S_{10}(T_2)$
$$\tau=4\,3\,2\,6\,7\,5\,10\,8\,9\,1,$$ we have
$$\alpha=4\,3\,2\,6\,7\,5,$$
with
$$M_1=4\quad M_2=6\quad M_3=7\quad M_4=10,$$
and
$$\beta=8\,9\,1.\\$$

\noindent The preceding considerations provide a characterization of the permutations in the two classes that end with the maximum symbol. Denote by $B_n$ the set of permutations in $S_n$ ending by $n$, by $B^{(1)}_n=B_n\cap Av(T_1)$, and by $B^{(2)}_n=B_n\cap Av(T_2)$. Theorem \ref{pincopallo} yields immediately the
following result:

\newtheorem{biie}[inizio]{Corollary}
\begin{biie}\label{forseavn}
The function $\psi_n:B_n\to S_{n-1}$ that maps a permutation
$\sigma$ into the permutation in $S_{n-1}$ obtained by deleting
the last symbol in $\sigma$ yields a bijection between \begin{itemize}
\item $B^{(1)}_n$ and
$S_{n-1}(321)$;
\item $B^{(2)}_n$ and
$S_{n-1}(312)$.
\end{itemize}
\end{biie}
\begin{flushright}
\vspace{-.5cm}$\diamond$
\end{flushright}

\section{Bijections with Dyck prefixes}

\noindent A \emph{Dyck prefix} is a lattice path in the integer
lattice $\mathbb{N}\times\mathbb{N}$ starting from the origin,
consisting of up steps $U=(1,1)$ and down steps $D=(1,-1)$, and
never passing below the $x$-axis. Obviously, a Dyck prefix can be
also seen as a word $W$ in the alphabet $\{U,D\}$ such that
every initial subword of $W$ contains at least as many symbols $U$
as symbols $D$.
\newline

\noindent It is well known (see e.g. \cite{vanl}) that the
number of Dyck prefixes of length $n$ is $\left({n\atop
\left\lfloor\frac{n}{2}\right\rfloor}\right)$.\newline

\noindent  A Dyck prefix ending at ground level is a \emph{Dyck
path}.\newline

\noindent We now define a map $\Phi:Av(T_1)\cup Av(T_2)\to \mathscr{P}$, where
$\mathscr{P}$ is the set of Dyck prefixes of even length.
First of all, associate the permutation $\sigma=1$ with the empty
path. Then, for every $n\geq 1$, associate a permutation
$\sigma\in S_{n+1}(T_1)\cup S_{n+1}(T_2)$ with a Dyck prefix of length $2n$, as
follows. Set $\sigma=M_1\,w_1\,M_2\,w_2\,\ldots\,M_k\,w_k$ as
above and let $l_i$ be the length of the word $w_i$. Now:

\begin{itemize}
\item if $w_k$ is empty, then $$\Phi(\sigma)=U^{M_1}D^{l_1+1}U^{M_2-M_1}D^{l_2+1}\cdots
U^{M_{k-1}-M_{k-2}}D^{l_{k-1}+1};$$
\item if $w_k=x_1\ldots x_{l_k}$ is not empty, then
$$\Phi(\sigma)=U^{M_1}D^{l_1+1}U^{M_2-M_1}D^{l_2+1}\cdots
U^{M_{k}-M_{k-1}}Q,$$ where $Q$ is the sequence $Q_1\ldots
Q_{l_k-1}$ of $l_k-1$ steps such that, for every
$j$, $Q_j$ is an up step if
$x_j=\max{\{x_j,x_{j+1},\ldots,x_{l_k}\}}$, a down step otherwise.
\end{itemize}

\noindent We point out that, in both cases, the last element of $\sigma$
is not processed. It is easy to check that the word $\Phi(\sigma)$
is a Dyck prefix of length $2n$.\newline

\begin{figure}[ht]
\begin{center}
\includegraphics[bb=77 592 488 719,width=.6\textwidth]{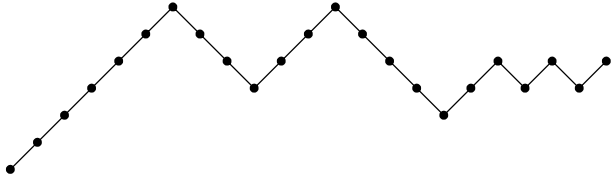} \caption{The
Dyck prefix $\Phi(6\ 1\ 2\ 9\ 3\ 4\ 5\ 11\ 12\ 7\ 10\
8)$.}\label{info}
\end{center}
\end{figure}

\noindent For example, consider the permutation in $S_{12}(T_1)$
$$\sigma=6\ 1\ 2\ 9\ 3\ 4\ 5\ 11\ 12\ 7\ 10\ 8.$$
We have $M_1=6$, $w_1=1\ 2$, $M_2=9$, $w_2=3\ 4\ 5$, $M_3=11$,
$w_3$ is empty, $M_4=12$, and $w_4=7\ 10\ 8$. The Dyck prefix
$\Phi(\sigma)$ is shown in Figure \ref{info}.

\noindent Note that, given a permutation $\sigma\in S_{n+1}(T_1)\cup S_{n+1}(T_2)$, $\sigma=M_1\,w_1\,\ldots\,M_k\,w_k$, the
position of the symbol $n+1$ (that is related to the existence and position of the $(n+1)$-th
up step in the associated Dyck prefix) plays an important role
in the definition of the map $\Phi$. For this reason, the $(n+1)$-th
up step in $\Phi(\sigma)$, if any, will be called the \emph{cut
step} of the path. Needless to say, $\Phi(\sigma)$ contains a cut
step if and only if it is not a Dyck path. In fact, if $w_k$ is empty, the path
$\Phi(\sigma)$ contains $M_{k-1}=n$ up steps, hence it is a Dyck
path. On the other hand, if $w_k$ is not empty, $\Phi(\sigma)$
contains at least $M_{k}=n+1$ up steps, therefore it does not end
at the ground level. The preceding considerations can be
summarized as follows:

\newtheorem{bleah}[inizio]{Proposition}
\begin{bleah}\label{hoepli}
Cosider a permutation $\sigma\in S_{n+1}(T_1)\cup S_{n+1}(T_2)$. Then,
$\Phi(\sigma)$ is a Dyck path if and only if $\sigma(n+1)=n+1$.
\end{bleah}
\begin{flushright}
\vspace{-.5cm}$\diamond$
\end{flushright}

\noindent Denote now by $\Phi_1$ and $\Phi_2$ the two restrictions of $\Phi$ to the sets $Av(T_1)$ and $Av(T_2)$. Our next goal is to prove that the two restrictions $\Phi_1$ and $\Phi_2$ are indeed bijections, by defining their inverses as follows. Both  $\Phi_1^{-1}$ and $\Phi_2^{-1}$ associate the
the empty path with the permutation $1$. Consider now a Dyck prefix
$P=U^{h_1}D^{s_1}U^{h_2}D^{s_2}\cdots U^{h_r}D^{s_r}$ of length
$2n$, $n\geq 1$. The permutation $\sigma=\Phi_1^{-1}(P)$ is defined
as follows:

\begin{itemize}
\item if $P$ is a Dyck path, namely, $h_1+\cdots +h_r=n=s_1+\cdots +s_r$,
set \begin{itemize}\item[]$\sigma(1)=h_1,$
 \item[]$\sigma(s_1+1)=h_1+h_2,$
 \item[]$\ldots$\item[] $\sigma(s_1+\cdots+s_{r-1}+1)=n,$
\item[]$\sigma(n+1)=n+1.$\end{itemize}
\noindent Then, place the remaining symbols in increasing order in
the unassigned positions.
\item if $P$ is not a Dyck path,  denote by $t$ the index of the ascending run $U^{h_t}$ containing the cut step. \begin{itemize} \item[-] if $t=1$, $P$ decomposes into
$$P=U^{n+1}Q,$$
where $Q$ is a lattice path. In this case, set
$$\sigma(1)=n+1;$$
\item[-] if $t>1$, $P$ decomposes into
$$P=U^{h_1}D^{s_1}\cdots
U^{h_t}D^{s_t}U^{n+1-h_1-\ldots-h_t}Q.$$ Set:\begin{itemize}
\item[] $\sigma(1)=h_1,$
\item[] $\sigma(s_1+1)=h_1+h_2,$
\item[] $\ldots$
\item[] $\sigma(s_1+\cdots+s_{t-1}+1)=n+1.$\end{itemize}
\end{itemize} In both cases, set $i=s_1+\cdots+s_{t-1}+1$ (or $i=1$ if $t=1$). Fill the unassigned positions less
than $i$ with the smallest remaining symbols placed in increasing
order. Then, for every $j=1,\ldots,n-i$, set either
\begin{itemize}
\item[] $\sigma(i+j)=\min [n+1]\setminus\{\sigma(1),\sigma(2),\ldots,\sigma(i+j-1)\} $ if
the $j$-th step of the path $Q$ is a down step, or
\item[] $\sigma(i+j)=\max [n+1]\setminus\{\sigma(1),\sigma(2),\ldots,\sigma(i+j-1)\}$ if
the $j$-th step of the path $Q$ is an up step.
\end{itemize}
Finally, $\sigma(n+1)$ equals the last unassigned symbol.
\end{itemize}

The permutation $\tau=\Phi_2^{-1}(P)$ can be defined
similarly:

\begin{itemize}
\item if $P$ is a Dyck path
set \begin{itemize}\item[]$\tau(1)=h_1,$
 \item[]$\tau(s_1+1)=h_1+h_2,$
 \item[]$\ldots$\item[] $\tau(s_1+\cdots+s_{r-1}+1)=n,$
\item[]$\tau(n+1)=n+1.$\end{itemize}
\noindent Then, scan the unassigned positions from left to right and fill them with the greatest unused symbol among those that are less then the closest preceding left-to-right maximum.
\item if $P$ is not a Dyck path, denote by $t$ the index of the ascending run $U^{h_t}$ containing the cut step.  \begin{itemize} \item[-] if $t=1$, $P$ decomposes into
$$P=U^{n+1}Q,$$
where $Q$ is a lattice path. In this case, set
$$\tau(1)=n+1;$$
\item[-] if $t>1$, $P$ decomposes into
$$P=U^{h_1}D^{s_1}\cdots
U^{h_t}D^{s_t}U^{n+1-h_1-\ldots-h_t}Q.$$ Set:\begin{itemize}
\item[] $\tau(1)=h_1,$
\item[] $\tau(s_1+1)=h_1+h_2,$
\item[] $\ldots$
\item[] $\tau(s_1+\cdots+s_{t-1}+1)=n+1.$\end{itemize}
\end{itemize} In both cases, set $i=s_1+\cdots+s_{t-1}+1$ (or $i=1$ if $t=1$). Then, scan the unassigned positions less then $i$ from left to right and fill them with the greatest unused symbol among those that are less than the closest preceding left-to-right maximum. Then, for every $j=1,\ldots,n-i$, set either
\begin{itemize}
\item[] $\tau(i+j)=\max \left([n+1]\setminus\{\tau(1),\tau(2),\ldots,\tau(i+j-1)\}\right)$ if
the $j$-th step of the path $Q$ is an up step, or
\item[] $\tau(i+j)=$ the second greatest element in the set \\ $[n+1]\setminus\{\tau(1),\tau(2),\ldots,\tau(i+j-1)\} $ if
the $j$-th step of the path $Q$ is a down step.
\end{itemize}
Finally, $\tau(n+1)$ equals the last unassigned symbol.
\end{itemize}

\noindent Theorem \ref{pincopallo} ensures that  $\sigma$ belongs to  $S_{n+1}(T_1)$, while  $\tau$ belongs to $S_{n+1}(T_2)$. Moreover,
it is easily seen that $\Phi_1^{-1}(\Phi_1(\sigma))=\sigma$ and $\Phi_2^{-1}(\Phi_2(\tau))=\tau$. As an
immediate consequence, we have:
\newtheorem{numero}[inizio]{Theorem}
\begin{numero}
The two maps $\Phi_1$ and $\Phi_2$ are bijections. Hence, the cardinality of both $S_{n+1}(T_1)$ and $S_{n+1}(T_2)$ is the central binomial
coefficient ${2n\choose n }$.\end{numero}
\begin{flushright}
\vspace{-.5cm}$\diamond$
\end{flushright}

\noindent We observe that the enumerative result contained in this theorem  can be also deduced from Theorem 4.6 in \cite{guib}.
\newline

\noindent In the following, we show that some properties of the permutations in $Av(T_1)$ and $Av(T_2)$ can be deduced
from certain features of the corresponding Dyck prefix.

\noindent First of all, a permutation $\sigma\in S_n$ is \emph{connected} if it has not a proper prefix of length $l<n$ that is a permutation in $S_l$.
Connected permutations appear in the literature also as \emph{irreducible permutations}.

\noindent On the other hand, recall that a \emph{return} of a Dyck prefix is a down
step ending on the $x$-axis. A Dyck prefix $P$ can be uniquely decomposed into $P=P'\,P''$, where $P'$ is a Dyck path and $P''$ is a floating Dyck prefix, namely, a Dyck prefix with no return (\emph{last return decomposition}). The last return decomposition of the path $\Phi(\sigma)$ gives information about the connected components of $\sigma$. More precisely, we have:

\newtheorem{ino}[inizio]{Proposition}
\begin{ino}\label{opportuno}
Let $\sigma$ be a permutation in $Av(T_1)\cup Av(T_2)$. The following are
equivalent: \begin{itemize}\item[a)] the Dyck prefix
$\Phi(\sigma)$ can be decomposed as $\Phi(\sigma)=P'\,P''$, where
$P'$ is a Dyck path of length $2l$, and $P''$ is a (possibly
empty) Dyck prefix,
\item[b)] $\sigma$ is the juxtaposition $\sigma'\,\sigma''$, where $\sigma'$ a
permutation of the set $\{1,\ldots,l\}$, and $\sigma''$ is a non-empty permutation.\end{itemize} In this case, letting $\tau$ be the permutation in $S_{l+1}(T_1)\cup S_{l+1}(T_2)$ obtained by
placing the symbol $l+1$ at the end of $\sigma'$, and
$\rho$ be the renormalization of
$\sigma''$, we have $P'=\Phi(\tau)$ and $P''=\Phi(\rho)$.
\end{ino}

\noindent \emph{Proof} We prove the assertion for permutations in $Av(T_1)$, the other case being analogous. Suppose that $\Phi(\sigma)$ can be
decomposed as follows:
$$\Phi(\sigma)=U^{h_1}D^{s_1}\cdots U^{h_p}D^{s_p}P'',$$
where $h_1+\cdots+h_p=s_1+\cdots+s_p$, and $P''$ is a Dyck prefix.
Set $l=s_1+\cdots+s_p$. By the definition of $\Phi_1^{-1}$, we have:
\begin{itemize}\item[]$\sigma(1)=h_1,$
\item[]$\sigma(s_1+1)=h_1+h_2,$
\item[]$\ldots$
\item[]$\sigma(s_1+\cdots+s_{p-1}+1)=l.$\end{itemize}

\noindent Now, we must fill the unassigned positions from $2$ to $l$ with
the $l-p$ smallest integers different from $h_1$, $h_1+h_2$,
$\ldots$, $l$. This implies that $\sigma(1)\,\ldots\,\sigma(l)$ is
a permutation of the set $\{1,\ldots,l\}$.

\noindent On the other hand, suppose $\sigma=\sigma'\,\sigma''$,
where $\sigma'$ is a permutation of the set $\{1,\ldots,l\}$, and
$\sigma''$ is a non-empty permutation. In this case,
$$\sigma=M_1\,w_1\,M_2\,w_2\,\ldots\,M_r\,w_r\,\sigma'',$$
where $M_r=l$. Note that the maximum symbol of $\sigma$ appears in
$\sigma''$. This implies that the portion of $\Phi(\sigma)$ that
corresponds to the entries in $\sigma'$ consists of $l$ up steps
and $l$ down steps, hence, it is a Dyck path.
\begin{flushright}
\vspace{-.5cm}$\diamond$
\end{flushright}

\begin{figure}[ht]
\begin{center}
\includegraphics[bb=56 572 379 663,width=.75\textwidth]{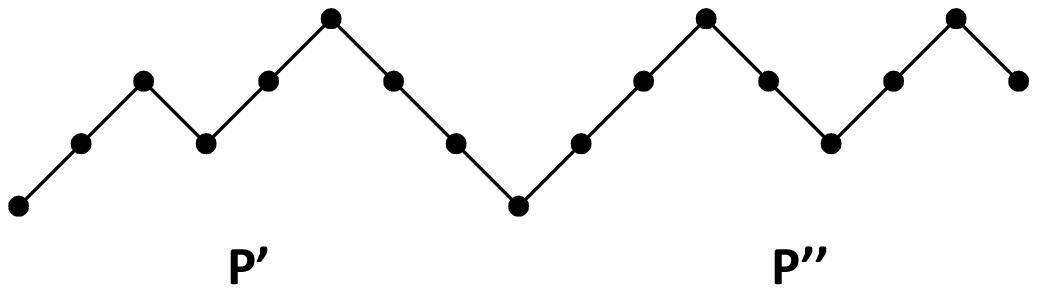} \caption{The
Dyck prefix $\Phi(2\,4\,1\,3\,7\,5\,9\,6\,8)$.}\label{tandc}
\end{center}
\end{figure}

\noindent For example, the path $P=P'\,P''$ in Figure \ref{tandc} corresponds
to the permutation $\sigma=\sigma'\,\sigma''$, where
$\sigma'=2\,4\,1\,3$ and $\sigma''=7\,5\,9\,6\,8$. Moreover, we
have $\tau=\Phi_1^{-1}(P')=2\,4\,1\,3\,5$ and
$\rho=\Phi_1^{-1}(P'')=3\,1\,5\,2\,4$.

\noindent Proposition \ref{opportuno} implies immediately the following result:

\newtheorem{lameta}[inizio]{Corollary}
\begin{lameta}
Connected permutations in $Av(T_1)$ (resp. $Av(T_2$)) are in bijection with floating Dyck prefixes. Hence, for every $n$, there are as many connected permutations in $S_n(T_1)$ (resp. $S_n(T_2)$) as non-connected permutations.
\end{lameta}

\noindent \emph{Proof} By Proposition \ref{opportuno}, we immediately deduce that connected permutations in $S_n(T_1)$ correspond bijectively to floating Dyck prefixes of length $2n-2$. These paths are in turn in bijection with Dyck prefixes of length $2n-3$ (one simply erases the first up step). Hence, the number of connected permutations in $S_n(T_1)$ is
$${2n-3\choose n-2}=\frac{|S_n(T_1)|}{2}.$$
\begin{flushright}
\vspace{-.5cm}$\diamond$
\end{flushright}

\noindent In closing, we consider the class $Av(T_1)\cap Av(T_2)=Av(T_1\cup T_2)$. We have:

\newtheorem{moreprec}[inizio]{Theorem}
\begin{moreprec}
$|S_n(T_1\cup T_2)|=n\cdot 2^{n-2}$.
\end{moreprec}

\noindent \emph{Proof} First of all, observe that a connected permutation $\sigma=\alpha\,n\,\beta$ belongs to $S_n(T_1\cup T_2)$ if and
only if $\alpha$ is an arbitrary increasing sequence not containing the symbol $1$ (since it must avoid both $321$ and $312$ and it must be connected) and $\beta$ is non-empty and either decreasing or order isomorphic to  $j\,j-1\,\ldots\,1\,2$ (since it must avoid $123$, $132$, $213$ et $231$). Denote by $k$ the length of $\alpha$. Then, the number of connected permutations in $S_n(T_1\cup T_2)$ is
$$\left(\sum_{k=0}^{n-2}2{n-2 \choose k}\right)-1=2^{n-1}-1.$$
Now, consider a permutation $\tau$ in $S_n(T_1\cup T_2)$ and decompose it as $\tau=\tau'\,\tau''$, where $\tau''$ is its longest connected suffix. Since $\tau''$ contains the symbol $n$, then $\tau'$ must avoid $321$ and $312$.
We distinguish the following cases:
\begin{itemize}
\item $\tau'$ is empty, hence $\tau$ is  connected. We have $2^{n-1}-1$ permutations of this kind.
\item $\tau'$ is non-empty and $\tau''$ contains at least two elements. If $k$ denotes the length of $\tau'$, we have $2^{k-1}(2^{n-1-k}-1)$ permutations of this kind.
\item $\tau''=n$. In this case $\tau$ avoids $321$ and $312$. We have $2^{n-2}$ permutations of this kind.
\end{itemize}
This means that:
$$|S_n(T_1\cup T_2)|=2^{n-1}-1+\left(\sum_{k=1}^{n-2}2^{k-1}(2^{n-1-k}-1)\right)+2^{n-2}=n\cdot 2^{n-2}.$$
\begin{flushright}
\vspace{-.5cm}$\diamond$
\end{flushright}

\section{Some statistics over the classes $Av(T_1)$ and $Av(T_2)$}

\noindent The definition of the map $\Phi$ suggests that some permutation statistics can be
studied simultaneously on the two sets $Av(T_1)$ and $Av(T_2)$:

\newtheorem{bsbsbs}[inizio]{Proposition}
\begin{bsbsbs}\label{nonsobene}
The three statistics  ``first element'',
``position of maximum symbol'', and ``number of left-to-right maxima'' are equidistributed over the classes $Av(T_1)$ and $Av(T_2)$.
\end{bsbsbs}

\noindent \emph{Proof} Consider a Dyck prefix $P$ of length $2n-2$ and the two permutations $\sigma=\Phi_1^{-1}(P)$
 and $\tau=\Phi_2^{-1}(P)$. Then, we can easily deduce the following:
 \begin{itemize}
 \item denote by $q$ the length of the first ascending run in $P$, namely the first maximal sequence of up steps. If $q\geq n$, then $\sigma(1)=\tau(1)=n$. Otherwise, $\sigma(1)=\tau(1)=q$;
 \item the position of $n$ in  both $\sigma$ and $\tau$ equals the number of down steps preceding the cut step, plus one;
 \item the left-to-right maxima different from $n$ in both $\sigma$ and $\tau$ correspond bijectively to peaks preceding the cut step.
 \end{itemize}
\begin{flushright}
\vspace{-.5cm}$\diamond$
\end{flushright}

\noindent Now we study the joint distribution of the two statistics ``position of maximum symbol'' and ``number of left-to-right maxima'' over the set $Av(T_1)$ (bearing in mind that this joint distribution is the same over $Av(T_2)$). More precisely, we determine the following generating function:
$$J(x,y,w)=\sum_{n\geq 1}\sum_{\sigma\in S_n(T_1)} x^ny^{pos(\sigma)}w^{lmax(\sigma)},$$
where $lmax(\sigma)$ denotes the number of left-to-right maxima in $\sigma$, and $pos(\sigma)$ denotes the position of the maximum symbol in $\sigma$.\\

\noindent In the study of permutation statistics over the two considered classes we exploit the last return decomposition of a Dyck prefix described in the previous section. The next Proposition analyzes the behavior of the statistics $pos(\sigma)$ and $lmax(\sigma)$ with respect to this decomposition:

\newtheorem{soddisfatta}[inizio]{Proposition}
\begin{soddisfatta}\label{nonmesso}
Consider a non connected permutation $\sigma\in Av(T_1)$. Consider the last return decomposition of $\Phi_1(\sigma)$
$$\Phi_1(\sigma)=P'\,P'',$$
where $P'$ is a non empty Dyck path and $P''$ is a floating Dyck prefix.
Set $\tau=\Phi_1^{-1}(P')$ and $\rho=\Phi_1^{-1}(P'')$. Then:
$$lmax(\sigma)=lmax(\tau)+lmax(\rho)-1,$$
$$pos(\sigma)=|\tau|+pos(\rho)-1.$$
\end{soddisfatta}

\noindent \emph{Proof} Proposition \ref{opportuno} implies that in
this case $\sigma=\sigma'\,\sigma'',$ where $\sigma'$ is obtained
from $\tau$ by deleting its last entry (which is a left-to-right
maximum), while $\sigma''$ is order isomorphic to $\rho$. For this reason, $lmax(\tau)=lmax(\sigma')+1$ and
$lmax(\rho)=lmax(\sigma'')$. Since the symbols appearing in
$\sigma''$ are greater than those appearing in $\sigma'$, we get
the first assertion. The second assertion is straightforward.
\begin{flushright}
\vspace{-.5cm}$\diamond$
\end{flushright}

\noindent For example, consider the path $P$ in Figure \ref{tandc}
and the permutation $\sigma=\Phi_1^{-1}(P)=2\,4\,1\,3\,7\,5\,9\,6\,8$. In this case,
$\tau=2\,4\,1\,3\,5$ and $\rho=3\,1\,5\,2\,4$, and
$$lmax(\tau)+lmax(\rho)-1=4=lmax(\sigma),$$
$$|\tau|+pos(\rho)-1=7=pos(\sigma).\vspace{.3cm}$$

\noindent The above result suggests to determine the joint
distribution of the two considered statistics over
 the set $B^{(1)}_n$ of permutations in $S_n(T_1)$ ending with the
maximum symbol, hence corresponding to Dyck paths, and over the set $C_n$ of connected permutations in $S_n(T_1)$, hence corresponding to floating Dyck prefixes. \newline

\noindent Set
$$B(x,y,w)=\sum_{n\geq 1}\sum_{\sigma\in B^{(1)}_n} x^ny^{pos(\sigma)}w^{lmax(\sigma)}.$$

\noindent As shown in the proof of Proposition \ref{nonsobene}, given a permutation $\sigma\in B^{(1)}_n$, we have $pos(\sigma)=|\sigma|=n$. Moreover, the number
of left-to-right maxima in $\sigma$ equals the number of peaks in $\Phi_1(\sigma)$, plus one. Hence, if we denote by $N(x,z)$ the Narayana function, namely,
the generating function of Dyck paths according to semi-length and number of peaks, then
$$B(x,y,w)=xywN(xy,w).$$
Exploiting the well-known expression of the Narayana function
$$N(x,z)=1+\frac{1-x(1+z)-\sqrt{(1-x(1+z))^2-4x^2z}}{2x}$$
(for more detailed information see, e.g., \cite{deu}), we get:
\begin{equation}\label{exprnara}
B(x,y,w)=w\frac{1+xy(1-w)-\sqrt{(1-xy(1+w))^2-4x^2y^2w}}{2}
 \end{equation}

\noindent Let now
$$C(x,y,w)=\sum_{n\geq 2}\sum_{\sigma\in
C_n}x^ny^{pos(\sigma)}w^{lmax(\sigma)}$$ be
the generating function of the joint distribution of $pos$ and $lmax$ over $C$. Note
that the summation above does not include the case $n=1$, since
the image under $\Phi$ of the unique permutation of length $1$ is the empty path, that is considered as
a Dyck path.

\noindent Proposition \ref{nonmesso}  yields the following functional equation
involving the generating functions $J(x,y,w)$, $B(x,y,w)$, and
$C(x,y,w)$:
\begin{equation} J(x,y,w)=B(x,y,w)+\frac{B(x,y,w)C(x,y,w)}{xyw}.\label{alfa}\end{equation}

\noindent Finally, we express the generating function $C(x,y,w)$ in
terms of $J(x,y,w)$ and $B(x,y,w)$. To this aim, we describe a
relation between the set of floating Dyck prefixes of length $2n$
and the set of all Dyck prefixes of length $2n-2$.

\noindent  Given a floating Dyck prefix $Q$, the lattice path
obtained from $Q$ by dropping its first and last step is a Dyck
prefix. On the other hand, given any Dyck prefix $P$, we can prepend to $P$ an up step and append either an up
or a down step, hence obtaining two Dyck prefixes $P_U$ and $P_D$,
respectively. The prefix $P_U$ is always floating, while $P_D$ is
floating if and only if the prefix $P$ is not a Dyck path.\newline

\noindent Denote now by $\sigma$ the permutation in $Av(T_1)$ corresponding to
a given Dyck prefix $P$ and suppose that $pos(\sigma)=h$ and $lmax(\sigma)=k$.

\begin{itemize}
\item if $P$ is floating, then both $P_U$ and $P_D$ are floating.
The definition of $P_U$ and $P_D$ implies that the cut steps in
both $P_U$ and $P_D$ correspond to the cut step in $P$. Set
$\sigma_U=\Phi_1^{-1}(P_U)$ and $\sigma_D=\Phi_1^{-1}(P_D)$. We have:
$$pos(\sigma_U)=h= pos(\sigma_D)$$
$$lmax(\sigma_U)=k=lmax(\sigma_D);$$
\item if $P$ is a Dyck path, only the path $P_U$ is floating. The
cut step in $P_U$ is obviously the last one. We have:
$$pos(\sigma_U)=h,$$
$$lmax(\sigma_U)=k.$$
\end{itemize}

\begin{figure}[ht]
\begin{center}
\includegraphics[bb=92 466 492 767,width=.6\textwidth]{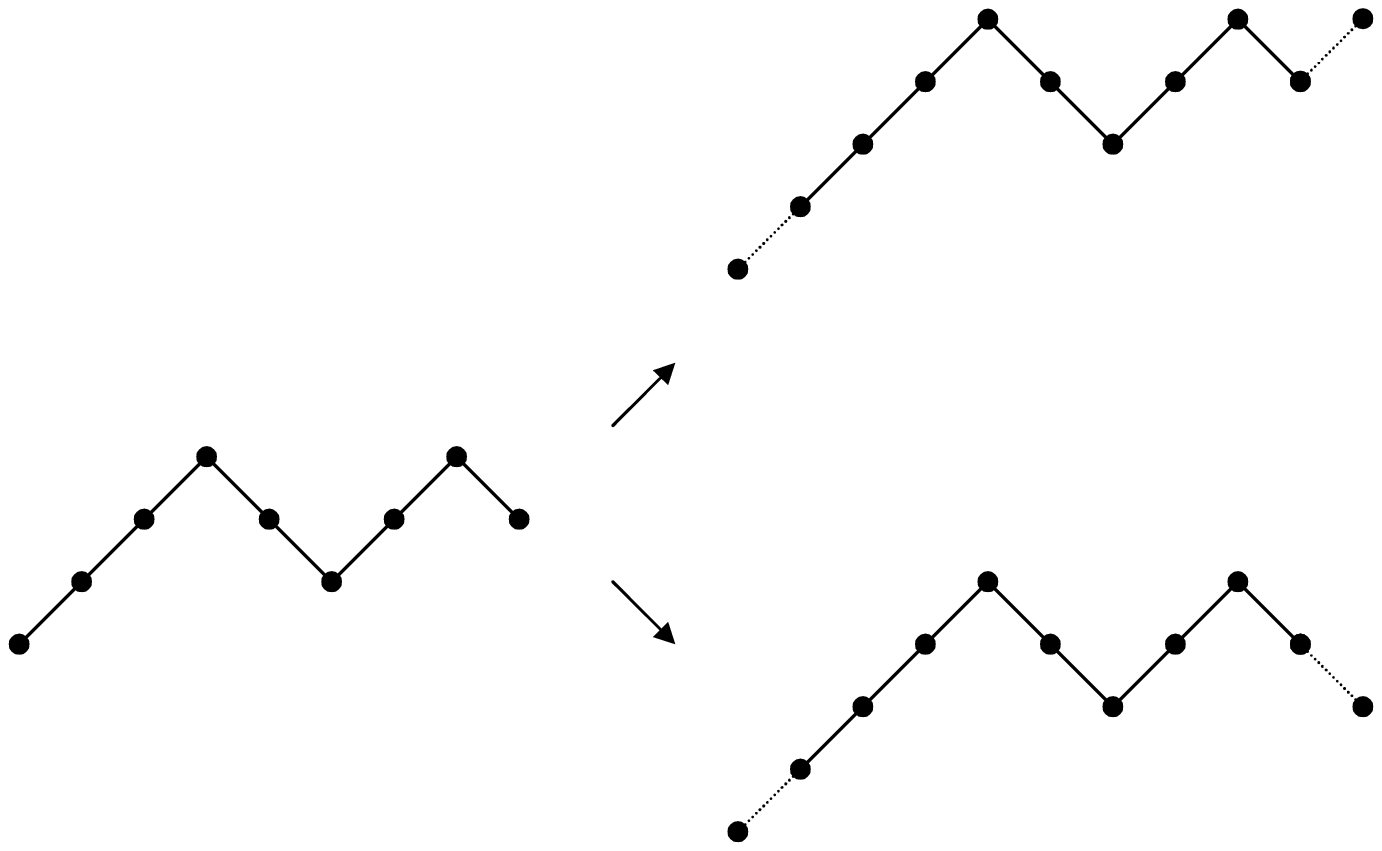}\caption{}\label{ehsi}
\end{center}
\end{figure}

\noindent For example, consider the permutations
$\sigma=\Phi_1^{-1}(P)=3\,1\,5\,2\,4$,
$\sigma_U=\Phi_1^{-1}(P_U)=4\,1\,6\,2\,5\,3$, and
$\sigma_D=\Phi_1^{-1}(P_D)=4\,1\,6\,2\,3\,5$, where $P$, $P_U$, and
$P_D$ are the Dyck prefixes in Figure \ref{ehsi}. We have
$pos(\sigma)=pos(\sigma_U)=pos(\sigma_D)=3$, and $lmax(\sigma)=lmax(\sigma_U)=lmax(\sigma_D)=2$.\\

\noindent Then we have:
\begin{equation}C(x,y,w)=2xJ(x,y,w)-xB(x,y,w).\label{beta}\end{equation}

\noindent Now, exploiting Identities (\ref{alfa}) and
(\ref{beta}), we get the following expression of $J(x,y,w)$ in terms
of $B(x,y,w)$:
\newtheorem{plinplin}[inizio]{Theorem}
\begin{plinplin} We have:
\begin{equation}J(x,y,w)=\frac{B(x,y,
w)(B(x,y,w)-yw)}{2B(x,y,w)-yw}.\label{taggala}\end{equation}
\end{plinplin}
\begin{flushright}
\vspace{-.5cm}$\diamond$
\end{flushright}
\noindent An explicit expression for $J(x,y,w)$ can be obtained by
combining Identities (\ref{exprnara}) and (\ref{taggala}).\newline

\noindent Let's now turn our attention to the statistic ``first element''. Given a permutation $\sigma$, we define $head(\sigma)=\sigma(1)$. We
determine the generating function
$$H(x,y)=\sum_{n\geq 1}\sum_{\sigma\in S_n(T_1)} x^ny^{head(\sigma)}=\sum_{n,k\geq 1}h_{n,k}x^ny^k,$$
where $h_{n,k}$ denotes the number of permutations $\sigma\in S_n(T_1)$ such that $head(\sigma)=k$ (by Proposition \ref{nonsobene}, this is also the generating function of the same distribution on $S_n(T_2)$).
First of all, given a permutation $\sigma\in S_n(T_1)$, if
$\sigma(1)=k$, then the Dyck prefix $\Phi(\sigma)$ starts with
\begin{itemize}
\item $k$ up steps followed by a down step, if $k<n$;
\item $n$ up steps, if $k=n$.
\end{itemize}

\noindent Hence, in the case $k<n$, if we delete from $\Phi_1(\sigma)$ the first
peak we obtain a Dyck prefix whose first ascending run contains at
least $k-1$ up steps. It is easy to see that this gives a
bijection between the set of Dyck prefixes of length $2n-2$
starting with $U^kD$ and the set of Dyck prefixes of length $2n-4$
starting with $U^t$, $t\geq k-1$. These arguments imply
that, if $n\geq 2$ and $k<n$:
\begin{equation}h_{n,k}=\sum_{j=k-1}^{n-1}h_{n-1,j}.\label{jkjk}\end{equation}
For $k>1$, this is equivalent to:
\begin{equation}h_{n,k}=h_{n,k-1}-h_{n-1,k-2},\label{equino}\end{equation}
with the convention $h_{s,0}=0$ for every integer $s$.\\

\noindent The special cases $k=1$ and $k=n$ can be easily handled as follows. First of all, the permutations
in $S_{n}(T_1)$ such that $\sigma(1)=1$ correspond to the Dyck prefixes of length $2n-2$ starting with $UD$, which are in one-to-one correspondence with Dyck prefixes
of length $2n-4$. Hence, $$h_{n,1}={2n-4\choose n-2}.$$

\noindent On the other hand, we observe that, given a Dyck prefix of length $2n-2$ starting with $U^n$, we can change
 the $n-$th up step into a down step, obtaining a new lattice path that is still a Dyck
 prefix. This implies that there are as many Dyck prefixes of
length $2n-2$ starting with $U^n$ as those starting with
$U^{n-1}D$, namely, $h_{n,n-1}=h_{n,n}$. Recall that permutations $\sigma\in
S_n(T_1)$ such that $\sigma(1)=n$ are in bijection with
permutations in $S_{n-1}(213,231)$. It is well known that the number of such permutations is $2^{n-2}$ (see \cite{simisch}). Hence,
$$h_{n,n-1}=h_{n,n}=2^{n-2}.$$

\newtheorem{incanato}[inizio]{Theorem}
\begin{incanato}
We have:
$$H(x,y)=\frac{xy\left[(xy-1)^2(1-y)\sqrt{1-4x}+x(1-2xy)\right]}{(1-y+xy^2)(1-2xy)\sqrt{1-4x}}.$$
\end{incanato}
\noindent \emph{Proof} Formula \ref{equino} gives a recurrence for the integers $h_{n,k}$ for every $n\geq 3$ and $2\leq k\leq n-1$. This fact suggests to
consider first the generating function
$$G(x,y)=\sum_{n\geq 2}\sum_{k=1}^{n-1}h_{n,k}x^ny^k.$$
Formula \ref{jkjk} yields:
$$G(x,y)=\sum_{n\geq 3}\sum_{k=2}^{n-1}h_{n,k-1}x^ny^k-\sum_{n\geq 3}\sum_{k=2}^{n-1}h_{n-1,k-2}x^ny^k+\sum_{n\geq 2}h_{n,1}x^ny=$$
$$=y\left(G(x,y)-x^2y-\sum_{n\geq 3}h_{n,n-1}x^ny^{n-1}\right)-xy^2\left(G(x,y)-\sum_{n\geq 2}h_{n,n-1}x^ny^{n-1}\right)+$$
$$+\sum_{n\geq 2}h_{n,1}x^ny.$$
The previous considerations allow us to deduce:
$$(1-y+xy^2)G(x,y)=\frac{x^3y^3-x^2y^2}{1-2xy}+\frac{x^2y}{\sqrt{1-4x}}.$$
Now, $H(x,y)$ can be obtained from $G(x,y)$ as follows:
$$H(x,y)=G(x,y)+xy+\sum_{n\geq 2} 2^{n-2}x^ny^n=G(x,y)+\frac{xy-x^2y^2}{1-2xy}.$$
\begin{flushright}
\vspace{-.5cm}$\diamond$
\end{flushright}

\section{Other statistics over $Av(T_1)$ and $Av(T_2)$}

\noindent This section is devoted to the study of some permutation statistics that are not equidistributed over the two classes. In both cases, we will translate occurrences of permutation statistics into configurations of the corresponding path. \\

\noindent First of all, we recall that a permutation $\sigma$ has an
\emph{ascent} at position $i$ whenever $\sigma(i)>\sigma(i+1)$,
and denote by $asc(\sigma)$ the number of ascents of $\sigma$.

\subsection{The class $Av(T_1)$}

\noindent We consider the generating function of the distribution of ascents over $Av(T_1)$:
$$F(x,y)=\sum_{n\geq 1}\sum_{\sigma\in S_n(T_1)} x^ny^{asc(\sigma)}.$$\newline

\noindent The ascents of $\sigma\in Av(T_1)$ can be recovered from the Dyck
prefix $\Phi_1(\sigma)$ as follows:

\newtheorem{collant}[inizio]{Proposition}
\begin{collant}
The number of ascents of a permutation $\sigma\in Av(T_1)$ is the
the sum of:
\begin{itemize}
\item the number of valleys and the number of triple descents (i.e.
occurrences of $DDD$) preceding the cut step in $\Phi_1(\sigma)$ (if $\Phi_1(\sigma)$ is
a Dyck path, its final down step counts as a valley) and
\item the number of down steps
following the cut step in $\Phi_1(\sigma)$.
\end{itemize}
\end{collant}

\noindent \emph{Proof} Decompose $\sigma$ as
$$\sigma=M_1\,w_1\,M_2\,w_2\,\ldots\,M_k\,w_k,$$ where $M_1,\ldots,M_k$ are the left-to-right maxima of $\sigma$.
Theorem \ref{pincopallo} implies that an ascent can occur in
$\sigma$ only in one of the following positions:
\begin{itemize}
\item between two consecutive symbols in $w_i$, $i\leq
k-1$. These two symbols correspond to two consecutive down steps
in $\Phi_1(\sigma)$ coming before the cut step. By the definition of
$\Phi_1$ these two down steps are necessarily preceded by a further
down step;
\item before every left-to-right
maximum $M_i$, except for the first one. These positions
correspond exactly to the valleys of $\Phi_1(\sigma)$ coming before
the cut step. In the special case when $\sigma$ ends with its
maximum symbol, the final ascent of $\sigma$ corresponds to the
final down step of the Dyck path $\Phi_1(\sigma)$;
\item in $w_k$, every time that the minimum unassigned symbol is
chosen. These ascents are of course in bijection with the down
steps following the cut step.
\end{itemize}
\begin{flushright}
\vspace{-.5cm}$\diamond$
\end{flushright}

\noindent Also in this case, we study the distribution of ascents on the set $B^{(1)}_n$ of permutations in $Av(T_1)$ such that $\Phi_1(\sigma)$ is a Dyck path, and on the set $C_n$ of connected permutations in $Av(T_1)$ that correspond to floating Dyck prefixes. Afterwards, we study the behavior of the ascent distribution with respect to the last return decomposition of the corresponding path. Arguments similar to those used in the proof of Proposition \ref{nonmesso} lead to:

\newtheorem{amal}[inizio]{Proposition}
\begin{amal}\label{deev}
Consider a permutation $\sigma\in Av(T_1)$. Suppose that
$\Phi_1(\sigma)$ can be decomposed into
$$\Phi_1(\sigma)=P'\,P'',$$
where $P'$ is a non empty Dyck path and $P''$ is any Dyck prefix.
Set $\tau=\Phi_1^{-1}(P')$ and $\rho=\Phi_1^{-1}(P'')$. Then:
$$asc(\sigma)=asc(\tau)+asc(\rho).$$
\end{amal}
\begin{flushright}
\vspace{-.5cm}$\diamond$
\end{flushright}

\noindent Consider the generating function of the ascent distribution over the set $B^{(1)}_n$:
$$E(x,y)=\sum_{n\geq 1}\sum_{\sigma\in B^{(1)}_n} x^ny^{asc(\sigma)}.$$
 In the paper \cite{bbs}, the authors determined the
generating function $A(x,y,z)$ of the joint distribution of
valleys and triple descents over the set of Dyck paths, namely,
\begin{eqnarray}
A(x,y,z) & = & \sum_{n\geq 0}\sum_{P\in \mathcal{P}_n}x^ny^{v(P)}z^{td(P)}=\nonumber\\ & = & \frac{1}{2xy(xyz-z-xy)}\left(-1+xy+2x^2y\right.
\nonumber\\
& &-2x^2y^2+xz-2xyz-2x^2yz+2x^2y^2z\\
& &
\left.+\sqrt{1-2xy-4x^2y+x^2y^2-2xz+2x^2yz+x^2z^2}\right)\nonumber
\end{eqnarray}

\noindent where $\mathcal{P}_n$ is the set of Dyck paths of semilenght $n$,
$v(P)$ denotes the number of valleys of the path $P$, and $td(P)$
is the number of occurrences of $DDD$ in $P$.\newline

\noindent We infer:

\begin{equation}E(x,y)=xy(A(x,y,y)-1)+x.\label{wewewei}\end{equation}
The last summand in Formula (\ref{wewewei}) takes into account the permutation $1$. This implies that:

\newtheorem{teorem}[inizio]{Proposition}
\begin{teorem}
We have:\begin{equation}
E(x,y)=\frac{\sqrt{1-4xy+4x^2y(y-1)}-1}{2y(x(y-1)-1)}\label{dentro}
\end{equation}
\end{teorem}
\begin{flushright}
\vspace{-.5cm}$\diamond$
\end{flushright}

\noindent  Consider now the generating function of the ascent distribution over the set $C_n$ of connected permutations in $Av(T_1)$
$$V(x,y)=\sum_{n\geq 2}\sum_{\sigma\in
C_n}x^ny^{asc(\sigma)}.$$

\noindent Proposition \ref{deev} yields the following functional equation
involving the generating functions $F(x,y)$, $E(x,y)$, and
$V(x,y)$:
\begin{equation} F(x,y)=E(x,y)+\frac{E(x,y)V(x,y)}{x}.\label{alfauno}\end{equation}

\noindent Finally, we express the generating function $V(x,y)$ in
terms of $F(x,y)$ and $E(x,y)$. Given any Dyck prefix $P$, we can obtain two Dyck prefixes $P_U$ and $P_D$ by prepending to $P$ an up step and appending either an up or a down step, as explained in the previous section and shown in Figure \ref{ehsi}.
In this case, we have:

\begin{itemize}
\item if $P$ is floating, then both $P_U$ and $P_D$ are floating, and
$$asc(\sigma_U)=h,\qquad asc(\sigma_D)=h+1,$$
\item if $P$ is a Dyck path, only the path $P_U$ is floating, and
$$asc(\sigma_U)=h.$$
\end{itemize}

\noindent Then, we have
\begin{equation}V(x,y)=(x+xy)F(x,y)-xyE(x,y).\label{betauno}\end{equation}

\noindent Now, exploiting Identities (\ref{alfauno}) and
(\ref{betauno}), we get the following expression of $F(x,y)$ in terms
of $E(x,y)$:
\newtheorem{plinpluc}[inizio]{Theorem}
\begin{plinpluc} We have:
\begin{equation}F(x,y)=\frac{E(x,y)(1-yE(x,y))}{1-E(x,y)-yE(x,y)}.\label{taggaladai}\end{equation}
\end{plinpluc}
\begin{flushright}
\vspace{-.5cm}$\diamond$
\end{flushright}
\noindent An explicit expression for $F(x,y)$ can be obtained by
combining Identities (\ref{dentro}) and (\ref{taggaladai}).\newline

\noindent In the remaining of this subsection, we characterize the permutations in $Av(T_1)$ according to the height of the last point of the path $\Phi_1(\sigma)$. Proposition \ref{hoepli} characterizes permutations in
$Av(T_1)$ whose associated prefix ends at the ground level. Now we
characterize permutations $\sigma \in S_n(T_1)$ whose corresponding
path $\Phi_1(\sigma)$ ends at $(2n-2,2h)$, $h>0$.

\noindent First of all, it is well known that the number of Dyck
prefixes of length $2n-2$ ending at $(2n-2,2h)$ is
\begin{equation}{2n-3\choose n-1-h}-{2n-3\choose
n-3-h},\label{finisci}\end{equation} (see \cite{vanl} and seq. A039599 in
\cite{oeis}).

\newtheorem{sequenza}[inizio]{Theorem}
\begin{sequenza}\label{longest}
Let $\sigma$ be a permutation in $Av(T_1)$ not ending with the
maximum symbol. If the $y$-coordinate of the last point of the
Dyck prefix $\Phi_1(\sigma)$ is $2k-2$, then the longest decreasing subsequence of $\sigma$ has cardinality $k$.
\end{sequenza}

\noindent \emph{Proof} Recall that every permutation $\sigma\in S_n(T_1)$ can be
decomposed as follows:
$$\sigma=\alpha\,n\,\beta,$$
where $\alpha$ avoids $321$ and $\beta$ avoids $213$ and
$231$. Moreover, $\beta=x_1\,x_2\,\ldots\,x_j$ is such that the
integer $x_i$ is either the minimum or the maximum of the set
$\{x_i,x_{i+1},\ldots,x_j\}$. Denote by $x_{i_1},\ldots,x_{i_q}$
the subsequence of $\beta$ consisting of the integers $x_i$ ($1\leq
i\leq j-1$) such that $x_i$ is the maximum of the set
$\{x_i,x_{i+1},\ldots,x_j\}$. By definition of the bijection
$\Phi$, it is immediately seen that the $y$-coordinate of the last
point of $\Phi_1(\sigma)$ is
$$j+1+q-(j-1-q)=2q+2.$$
It is easy to check that the sequence
$$n\,x_{i_1}\,\ldots\,x_{i_q}\,x_j,$$ of length $k=q+2$, is the longest decreasing
subsequence in $\sigma$. This ends the proof.
\begin{flushright}
\vspace{-.5cm}$\diamond$
\end{flushright}

\noindent The preceding result allows us to characterize the set
of Dyck prefixes of length $2n-2$ corresponding via $\Phi_1$ to
permutations in $S_n(T_1)$ that avoid also the pattern
$k\,\,k-1\,\ldots\,2\,1$:

\newtheorem{gigigi}[inizio]{Theorem}
\begin{gigigi}
We have:
$$|S_n(T_1,k\,\,k-1\,\ldots\,2\,1)|={2n-2 \choose n-1}-{2n-2\choose n-k}$$
\end{gigigi}

\noindent \emph{Proof} The preceding results yield immediately:
$$|S_n(T_1,k\,\,k-1\,\ldots\,2\,1)|=\sum_{i=0}^{k-2}{2n-3 \choose n-1-i}-{2n-3\choose n-3-i}=$$$$={2n-2 \choose n-1}-{2n-2\choose n-k}.$$

\begin{flushright}
\vspace{-.5cm}$\diamond$
\end{flushright}

\noindent In particular, consider the case $k=3$. Of course, we have $S_n(T_1,321)=S_n(321)$. The set of Dyck prefixes of length $2n-2$
corresponding via $\Phi_1$ to permutations in
$S_n(321)$ can be partitioned into two subsets:
\begin{itemize}
\item[a)] the set of Dyck paths;
\item[b)] the set of Dyck prefixes ending at $(2n-2,2)$.
\end{itemize}
\noindent It is well known that the set $S_n(321)$ is enumerated by $n$-th Catalan number. Many bijections between permutations in
$S_n(321)$ and Dyck paths of semilength $n$ appear in the literature, notably the bijection defined by
Krattenthaler in \cite{kratt}. If $\sigma$ is a permutation in $S_n(321)$, the relation between the Dyck prefix $\Phi_1(\sigma)$ and the Dyck path
$K(\sigma)$ associated with $\sigma$ by Krattenthaler's bijection can be described as follows:
\begin{itemize}
\item if $\Phi_1(\sigma)$ is a Dyck path, $K(\sigma)=\Phi_1(\sigma)UD$;
\item if the last point of $\Phi_1(\sigma)$ has coordinates $(2n-2,2)$, $K(\sigma)=\Phi_1(\sigma)DD$.
\end{itemize}

\section{The class $Av(T_2)$}

\noindent In this last section, we study the generating function of the ascent distribution over $Av(T_2)$
$$M(x,y)=\sum_{n\geq 1}\sum_{\sigma\in S_n(T_2)} x^ny^{asc(\sigma)}.$$

\newtheorem{kio}[inizio]{Proposition}
\begin{kio}
The number of ascents in a permutation $\sigma\in Av(T_2)$ is the
number of peaks in the Dyck prefix $\Phi_2(\sigma)$.
\end{kio}

\noindent \emph{Proof} Decompose $\sigma$ as
$$\sigma=M_1\,w_1\,M_2\,w_2\,\ldots\,M_k\,w_k,$$ where $M_1,\ldots,M_k$ are the left-to-right maxima of $\sigma$.
Theorem \ref{pincopallo} implies that an ascent can occur in
$\sigma$ only in one of the following positions:
\begin{itemize}
\item before every left-to-right
maximum $M_i$, except for the first one. These positions
correspond exactly to the peaks of $\Phi_2(\sigma)$ coming before
the cut step.
\item in $w_k$, every time that the second greatest unassigned symbol $x_j$ is
chosen either immediately after the cut step or immediately after the choice of the maximum unassigned element $x_{j-1}$. This means that there is an element $p$ such that $x_j<p<n$ or $x_j<p<x_{j-1}$, respectively. In both cases, when $p$ will be placed, it will give rise to an ascent in $\sigma$.
These ascents are easily seen to be in bijection with peaks following (or involving) the cut step.
\end{itemize}
\begin{flushright}
\vspace{-.5cm}$\diamond$
\end{flushright}

\noindent Hence, we study the distribution of peaks on the set $\mathscr{P}$ of Dyck prefixes, namely, the generating function
$$S(x,y)=\sum_{n\geq 0}\sum_{P\in \mathcal{P}_n}x^ny^{peak(P)},$$
where $peak(P)$ denotes the number of peaks in the prefix $P$.

\noindent Denote by $R(x,y)$ and $\hat{R}(x,y)$ the generating functions of the same distribution on the set of floating Dyck prefixes and on the set of Dyck prefixes ending with $U$, respectively.
We have:

\newtheorem{venticinque}[inizio]{Proposition}
\begin{venticinque}
The two generating functions $R(x,y)$ and $\hat{R}(x,y)$ coincide.
\end{venticinque}

\noindent \emph{Proof} Consider a floating Dyck prefix $P$. We can obviously write $P=U\,P'$, where $P'$ is still a Dyck prefix. Consider now the path $Q=P'\,U$, ending with $U$. Then, the map $P\mapsto Q$ is a size-preserving bijection between the set of floating Dyck prefixes and the set of Dyck prefixes ending by $U$, such that $peak(P)=peak(Q)$.
\begin{flushright}
\vspace{-.5cm}$\diamond$
\end{flushright}


\noindent Consider now a Dyck prefix $Q$ ending with $U$. Then, according to the last return decomposition, we can decompose $Q$ into $Q=P\,U\,P'\,U$, where $P$ is a Dyck path and $P'$ is any Dyck prefix. We deduce that:
\begin{equation}\label{ldtanto} \hat{R}(x,y)=xN(x,y)S(x,y)=R(x,y),
\end{equation}
where $N(x,y)$ is the Narayana generating function.

\noindent Afterwards, exploiting once again the last return decomposition, we have:

\begin{equation}
\label{ioiuy}S(x,y)=N(x,y)\left(R(x,y)+1\right).
\end{equation}

\noindent Then, combining Identities (\ref{ldtanto}) and (\ref{ioiuy}), we obtain

$$S(x,y)=\frac{N(x,y)}{1-xN(x,y)^2}.$$

\noindent Consider now the generating function $M(x,y)$ of the ascent distribution over $Av(T_2)$.
The definition of the map $\Phi_2$ yields immediately that $M(x,y)=xS(x,y)$. Hence:

\newtheorem{ottomin}[inizio]{Theorem}
\begin{ottomin}
We have:
$$M(x,y)=\frac{N(x,y)}{1-xN(x,y)^2}$$
\end{ottomin}
An expression for $M(x,y)$ can be found by replacing $N(x,y)$ by its explicit formula.


\end{document}